\documentclass{article}

\usepackage{amssymb, amsmath, amsthm, enumitem, multirow, cancel, physics, color, mathtools, xcolor, mathrsfs , hyperref,tkz-euclide, tikz, bigints}

\usepackage{mathtools}
\DeclarePairedDelimiter\ceil{\lceil}{\rceil}

\usepackage{stackengine}

\usepackage{arxiv}

\usepackage[utf8]{inputenc} 
\usepackage[T1]{fontenc}    
\usepackage{hyperref}       
\usepackage{url}            
\usepackage{booktabs}       
\usepackage{amsfonts}       
\usepackage{nicefrac}       
\usepackage{microtype}      
\usepackage{lipsum}

\title{A new Laplace-type fractional derivative}

\author{
Mostafa Rezapour \\
  Department of Mathematics\\
  Washington State University\\
  Pullman WA, 99163 \\
  \texttt{mostafa.rezapour@wsu.edu} \\
  \And
 Adebowale Sijuwade\\
  Department of Mathematics\\
  Washington State University\\
  Pullman WA, 99163 \\
  \texttt{adebowale.sijuwade@wsu.edu} \\
}

\begin{document}
\maketitle

\begin{abstract}

In this paper, we present a new derivative via the Laplace transform. The Laplace transform leads to a natural form of the fractional derivative which is equivalent to a Riemann-Liouville derivative with fixed terminal point. We first consider a representation which interacts well with periodic functions, examine some rudimentary properties and propose a generalization. The interest for this new approach arose from recent developments in fractional differential equations involving Caputo-type derivatives and applications in regularization problems.

\end{abstract}


\section{Introduction}

This section is intended to motivate the usage of fractional differentiation and integration. Let D denote the ordinary derivative $\frac{d}{dt}\footnotesize( \cdot \footnotesize)$. There are many applications of fractional calculus such as: modeling diffusion processes in Atananckovic \cite{Ref1}; heat transfer models in Atangana and Baleanu \cite{Ref3}; continuum mechanics in Mainardi \cite{Ref4}; evolution equations and renewal processes in Kochubei \cite{Ref5}; physics and engineering in Diethelm \cite{Ref8}; robotics, aerospace and biomedicine in Caponetto \cite{Ref13} and financial economics in Kiryakova \cite{Ref16}. Partial fractional differential equations and modelling techniques in hydrodynamics and stochastics can be found in Kilbas et al. \cite{Ref15}. An introductory treatment of fractional differential equations is provided in Miller et al. \cite{Ref18} and Podlubny \cite{Ref20}. One can find several different definitions involving fractional-order integral and derivative operators such as those seen in de Olivera \cite{Ref2} and Samko \cite{Ref19}.  In this article, we consider spaces where the functions are continuous or piecewise continuous.

\bigskip

Suppose that $f$ is locally integrable on $\mathbb{R}$. We define distributions by

\[ < f(t), \phi(t) > = \int_{\mathbb{R}} f(t) \phi(t) dt,  \tag{1} \label{eq: 1}
 \]

\noindent and denote the space of distributions by $\mathscr{D}$'. Consider the delta distribution in the usual sense as a continuous linear functional acting on the space of distributions $\mathscr{D}$ by

\[ <\delta(t-c), \phi(t) > = \phi(c) \hspace{1mm} \text{ for } \phi \in \mathscr{D}, c \in \mathbb{R}.  \]

The derivative of a distribution $\phi$ is given by 

\[ <\phi’(t), \psi(t)> = -<\phi(t), \psi’(t)>, \]  

\noindent where $\psi \in \mathscr{D}$. Setting $\phi(t)=\delta(t-c)$, we obtain $\delta’(c)$ on the right hand side. This definition can be easily generalized to higher order derivatives by 

\[ <D^{(n)} \delta(t-c) , \phi(t)> = (-1)^{(n)} (D^{(n)} \phi)(c). \]

\noindent The Fourier transform of a distribution is defined as follows:

\[ <  \mathscr{F} (\phi(t)) (\omega), \phi(\omega)> = < \phi(t), \mathscr{F} (\phi(\omega))>, \] 
\noindent where $\mathscr{F} (\phi(t))$ is the usual definition of the Fourier transform for real-valued functions. 

\bigskip 

Now consider the Schwarz class of smooth test functions $\mathscr{I}$, which decay with their derivatives at infinity. The space of continuous linear functionals on this space is denoted by  $\mathscr{I}$’, which is contained in the set of distributions. Consider the Laplace transform of a distribution $\phi(t)$ defined by

\[ F(s) = \mathscr{L}(\phi(t))= \mathscr{F} (\phi(t)e^{-\sigma t})(\mu), \tag{2} \label{eq: 2}
 \]

\noindent where $s = \sigma + i\mu$, $ \mu < 0$ and $\phi(t)e^{-\sigma t} \in \mathscr{I}'$. Let $f$ be a distribution supported on $(0, \infty)$ such that for $\sigma >0$, $f(t)e^{-\sigma t} \in \mathscr{I}'$. It follows that the Laplace transform of the Dirac delta function is given by 

\[ \mathscr{L}(\delta(t-c)) (s) = e^{-cs}, \] 

\noindent and

\[ \mathscr{L}(\phi’(t))(s) = s\mathscr{L}(\phi(t)) (s). \] 

A more rigorous treatment of distributional derivatives and generalized functions can be found in Schwarz \cite{Ref6}, Kanwal \cite{Ref7},  McBride \cite{Ref17} and Zemanian \cite{Ref30}. If $f:\mathbb{R} \to \mathbb{R}$ is of exponential order, then

\[ F(s)=\mathscr{L}f(t) = \int_{0}^{\infty} e^{-st} f(t) dt, \] 

\[  \mathscr{L}(D^{n} f) = s^{\alpha} F(s) - \sum\limits_{k=0}^{n-1} s^{k} (D^{n-k-1} f(t) )_{t=0} . \tag{3} \label{eq: 3} \]

The left and right Riemann Liouville derivatives of order $\alpha \in \mathbb{C}$, where $Re(\alpha)\geq 0$ on a finite real interval are defined by 

\[ (^{L}_{a} \mathbf{D} ^{\alpha}_{t})g (t) = \frac{1}{\Gamma(n-\alpha)} D^{n} \int_{a}^{t} g(s) (t-s)^{n-\alpha-1} ds, \tag{4} \label{eq: 4}  \]

when $t>a$ and

\[ (^{R}_{b} \mathbf{D} ^{\alpha}_{t})g (t) = \frac{1}{\Gamma(n-\alpha)} D^{n} \int_{t}^{b}  g(s) (t-s)^{n-\alpha-1} ds, \tag{5} \label{eq: 5}  \]

when $t < b$ that resemble the Cauchy formula for integration given by

\[ (D^{-n} f)(t) = \int_{a}^{t} \frac{f(s)(t-s)^{\alpha-1} ds}{\Gamma(\alpha)},  \]

where $t<b$, $n=\ceil{\alpha}$. From the identity

\[ D^n f(t) = \lim\limits_{h \to 0} h^{-n} \sum\limits_{k=0}^{n} (-1)^{k} {n \choose k} f(t-kh),  \]

we have the Grunwald-Letnikov definition for $\alpha <0$: 

\[  ^{GL}_a D_{t}^{-\alpha} f(t) = \lim_{\substack{h \to 0 \\ nh = t - a}}  h^{-\alpha} \sum\limits_{k=0}^{n} (-1)^{k} {\alpha \choose k} f(t-kh) \tag{6} \label{eq: 6},  \]

and it can be shown that if $f \in C[a,t]$,

\[ \lim_{\substack{h \to 0 \\ nh = t - a}} h^{-\alpha} \sum\limits_{k=0}^{n} (-1)^{k} {\alpha \choose k} f(t-kh)  \]

\[ = \frac{1}{\Gamma(\alpha)} \lim\limits_{n, k \to \infty}  \sum\limits_{k=0}^{n} \bigg( \frac{(-1)^k \Gamma(\alpha)}{k^{\alpha - 1}} {-\alpha \choose k} \bigg) \frac{t-a}{n} \bigg (\frac{k(t-a)}{n}\bigg)^{\alpha -1} f\bigg( t-\frac{k(t-a)}{n} \bigg) \]

\[ = \frac{1}{\Gamma(\alpha)} \int_{a}^{t} (t-s)^{\alpha - 1} f(s) ds, \]

\noindent which connects the Grunwald-Letnikov and Riemann-Liouville approaches. For $\alpha<0$ it is then enough to replace ${ \alpha \choose k }$ with  $(-1)^{-k} { -p \choose k }$ in the classical limit and proceed. The Riemann-Liouville fractional derivative is problematic in that 

\[  \lim\limits_{t \to a} {^{L}}_{a} \mathbf{D} ^{\alpha}_{t} = b_n,   \]

\noindent where $b_n$ for $n=0,1,2, ...$ are prescribed constants. Consider the $\alpha$-th Caputo fractional derivative of $f(t)$ defined by

\[ ^C_a D_{t}^{\alpha} f(t) = \frac{1}{\Gamma(\alpha - n)} \int_{a}^{t} \frac{f^{(n)} (s) ds}{(t-s)^{\alpha + 1 - n}},  (n-1 < \alpha < n). \tag{7} \label{eq: 7} \]

\noindent Integrating by parts, one can see that for $f \in C^{n+1} ([a,T])$, $0 \leq n-1 < \alpha < n$,

\[  ^C_a D_{t}^{\alpha} f(t) \to D^{n} f(t) \]

\noindent as $\alpha \to n$. The Caputo derivative unfortunately does not reduce to the classical derivative for $n-1 < \alpha \leq n$ since its behavior depends on the terminal point, $a$ as follows


\[ ^C_a D_{t}^{k} f(t) \to  \int_{a}^{t} \frac{f^{(n)} (s) ds}{\Gamma(1) } = D^{n-1} f(t) - D^{n-1} f(a). \]

The Caputo and Riemann-Liouville derivatives do not satisfy all of the expected classical properties. For instance, for $c \in \mathbb{R}$, $n-1 < \alpha \leq n$,

\[  ^C_a D_{t}^{\alpha} c = 0 \text{ while } ^{L} _a  \textbf{D}_{t}^{\alpha} c = \frac{ct^{-\alpha}}{\Gamma(\alpha - n)}. \] 



\noindent Fortunately, the left Riemann Liouville fractional derivative does satisfy a semigroup property under the following conditions. Suppose that $\alpha \neq \beta$. If $f$ is $N$ times differentiable where $N=\max{ (\ceil{\alpha}, \ceil{\beta}) }$ and $D^{k} f(a)=0$ for $k=0,1,2, … N$, we have that

\[  ^{L}_{a} \mathbf{D}^{\alpha}_{t} (^{L}_{a} \mathbf{D}^{\beta}_{t} f(t) ) = ^{L} _a  \textbf{D}_{t}^{\alpha + \beta},  \]

\noindent as in equation (2.127) of \cite{Ref20}. In the case of the Caputo fractional derivative and Riemann-Liouville fractional derivatives, if $D^{s} f (0) = 0$ for $s=n, n+1 , … m$, we have

\[ ^{C} _{a} D^{\alpha}_{t}  (^{C} _{a} D^{m}_{t} f(t) ) = ^{C}_{a} D^{\alpha + m}_{t} f(t),  \hspace{1mm}  m=0,1,2, … \]  

\noindent and 

\[ ^{L} _{a} \textbf{D}^{m}_{t}  (^{L} _{a} \textbf{D}^{\alpha}_{t} f(t) ) = ^{L} _{a} \textbf{D}^{\alpha + m}_{t} f(t), \hspace{1mm}  m=0,1,2, … \]  

\noindent For the Riemann-Liouville and Grunwald-Letnikov fractional derivatives, we also have the expected result

\[ D^{\alpha} (t-c)^{\gamma}= \frac{\Gamma(\gamma + 1)}{\Gamma(\gamma - \alpha + 1 )}(t-c)^{\gamma - \alpha}, \tag{8} \label{eq: 8} \]

\noindent The equivalence of the Riemann-Liouville and integral form of the Grunwald-Letnikov fractional derivative (\cite{Ref20}  pg. 63) allows us to verify this as follows. Since
 
\[ ^{GL}_a {D}_{t}^{\alpha} f(t) = \frac{1}{\Gamma(-\alpha)} \int_{a}^{t} (t-s)^{-\alpha - 1} (s-c)^{\gamma} ds,  \]

\noindent substituting $u = \frac{s-c}{t-c}$, we have:




\[ \int_{0}^{1}  \frac{ ((t-c)(1-u))^{-\alpha-1}    ((t-c)u)^{\gamma}  (t-c)}{\Gamma(-\alpha)}  du  \]

\[  = \frac{(t-c)^{\gamma - \alpha}}{\Gamma(-\alpha)}   \int_{0}^{1}  { (1-u)^{-\alpha-1}  u^{\gamma}}  du  \]

\[ =    (t-c)^{\gamma - \alpha} \bigg( \frac{\beta(\gamma + 1, -\alpha)} {\Gamma(-\alpha)} \bigg)  \]




\[ = \frac{\Gamma(\gamma + 1)}{\Gamma(\gamma - \alpha + 1 )}(t-c)^{\gamma - \alpha}.  \]



\noindent These definitions can be extended to higher dimensions as in Caputo and Fabrizio \cite{Ref22}. More details on the Caputo derivative can be found in Caputo \cite{Ref14}. Applications of these fractional derivatives in the context of gradient descent algorithms can be found in Wei et al \cite{Ref24} and Pu et al. \cite{Ref25}

\section{The $\ell_{1}$ derivative}
\label{sec:1}

In this paper, we seek a fractional derivative using a transformation approach. Suppose that $f^{(k)}(0)=0$ for $k=0,1,2 , ... , n-1$. A Laplace-transform based fractional derivative can be defined naturally by

\[ ^{LT} D^{\alpha}_t  f = \mathscr{L}^{-1} ( s^{\alpha}F(s) ). \tag{9} \label{eq: 9} \] 

\noindent In the case that $f= (t-c)^{\gamma}$, we have 

\[ ^{LT} D^{\alpha}_t  (t-c)^{\gamma} = \mathscr{L}^{-1} ( s^{\alpha}\mathscr{L}( (t-c)^{\gamma} ) ) \]


\[ = \Gamma(\gamma + 1) \mathscr{L}^{-1} ( e^{-sc}  s^{\alpha - \gamma - 1}  ) \] 

\[ =  \frac{\Gamma(\gamma + 1)} {\Gamma(\gamma - \alpha + 1) } \mathscr{L}^{-1}( e^{-sc}  \mathscr{L}(t^{\gamma - \alpha}) ) \]   

\[ = \frac{\Gamma(\gamma + 1)} {\Gamma(\gamma - \alpha + 1) }  (t-c)^{\gamma - \alpha} H(t-c), \] 

\noindent where $H(t)$ is the usual Heaviside function. For simplicity, let $f(t) = | t - c | $, then

%

%




\[  ^{LT} D^{\alpha}_t f(t) = \mathscr{L}^{-1}(s^{\alpha} \mathscr{L}( (c-t) + 2(t-c)H(t-c)) ), \] 

\noindent since

\[ \mathscr{L}^{-1}(s^{\alpha} \mathscr{L}( (c-t) + 2(t-c)H(t-c)) ) \]

\[ = \mathscr{L}^{-1} ( cs^{\alpha-1} +  s^{\alpha-2}(2e^{-cs} - 1)     ).   \]

\noindent When $\alpha = 1$, this reduces to $2H(t-c) - 1$, which agrees with the usual distributional derivative $\frac{t}{|t|}$. When $\alpha = 2$, 

\[ ^{LT} D^{\alpha}_t |t-c| = (cD^{*} - 1)\delta(t) + 2H(t-c), \] 

\noindent otherwise when $\alpha <1$, 

\[ ^{LT} D^{\alpha}_t |t-c| =  \frac{ct^{\alpha - 2}}{\Gamma(\alpha -1)} - \frac{ct^{\alpha - 3}(2H(t-3) - 1)}{\Gamma(\alpha -2)}. \] 

\noindent For other values of $\alpha$, the inverse Laplace transform does not exist, since by a standard choice of contour as in Arfken [34], for $p>0$,

\[  \frac{1}{2\pi i} \int_{\gamma- i \sqrt{R^2  - \gamma^2} }^{\gamma +  i \sqrt{R^2  - \gamma^2} }  s^{p}e^{st} ds = - \frac{1}{2\pi i} \int_{ \frac{3\pi}{2} - \sin^{-1}(\frac{\gamma}{R}) }^{\frac{\pi}{2} + \sin^{-1}(\frac{\gamma}{R})}  (Re^{i\theta})^pe^{Rte^{i \theta}} iRe^{i\theta} d\theta. \tag{10} \label{eq: 10} \]

\noindent After integration by parts and a standard application of the Euler's reflection formula, the above integral results in

\[ \int_{-\epsilon}^{0} (-se^{-i\pi})^{p} e^{st} \frac{ds}{2\pi i }  - \frac{1}{2\pi i} \int_{-\epsilon}^{0} (-se^{i\pi})^{p} e^{st} \frac{ds}{2\pi i }  + \int_{-\pi}^{\pi} (\epsilon e^{i\theta})^p i \epsilon e^{i\theta} \frac{ds}{2\pi i}   \]

\[ = \frac{\sin{p\pi}} {\pi} \int_{\epsilon}^{\infty} e^{-st} s^{p} ds + \frac{\epsilon^{p+1}\sin{\pi p} }{\pi(p+1)}    \] 

\[ = \frac{t^{-p-1} }{\Gamma(-p)},  \]  
\noindent as $R \to \infty$, $ \epsilon \to 0$. In the case of the Riemann-Liouville fractional derivative, we obtain

\[ \mathscr{L} (_{0} D^{\alpha}_{t} f(t) ) = s^{\alpha} \mathscr{L}f(t)  - \sum\limits_{k=0}^{n-1} s^{k} (D^{n-k-1} f(t) )_{t=0} \]

\[ =s^{\alpha} F(s) - \sum\limits_{k=0}^{n-1} D^{n-1-k} ( g_{n-\alpha} \star f) )_{t=0},  \] 

\noindent where $n = \ceil{\alpha}$. Thus, the Laplace transform readily generalizes the left Riemann-Liouville fractional derivative with a lower limit $a=0$, applying ($\ref{eq: 3}$). When $\alpha = \pm 1$,

\[  \mathscr{L}(D^{n} f(x) ) = \mathscr{L}^{-1} ( s^{n} \mathscr{L}(f) ), \] 

\noindent reduces to

\[ s\mathscr{L}\big( \int_{0}^{t} f(x) dx \big) =  \mathscr{L}(f(t)). \]


\noindent Now let 

\[ g_{c}(t) = \frac{t^{c - 1}H (t)}{\Gamma(c)}, \]

\noindent then

\[   D^{\alpha}_{t} f(t) = \mathscr{L}^{-1} ( s^{\alpha} \mathscr{L} f(t) )  \]

\[ = \mathscr{L}^{-1} (s^{n} \mathscr{L} (g_{n-\alpha} \star f)(t)) \] 
\[ = \mathscr{L}^{-1} (\mathscr{L} (g_{n-\alpha} ) (sF(s))) \]

\[  =  g_{n-\alpha}(t) \star Df(t), \tag{10} \label{eq: 10} \]

\noindent which agrees with the left Riemann-Liouville derivative for $a =0$.

We now consider a new Laplace derivative. Suppose that $\sigma >0$ and $f$ is a real-valued function such that $\{ f(t) e^{-\sigma t} \} \in \mathscr{I}'$. 

\bigskip

\noindent \textbf{Definition 2.1.} Define the $\ell_{1}$ fractional derivative by

\[ ^{\ell_{1}} (D f)(t) =  {\lim\limits_{\epsilon \to 0}} \mathscr{L}^{-1}\bigg(s \bigg( \mathscr{L}(f(t-\epsilon) )- f(t\epsilon) \bigg) \bigg), \tag{11} \label{eq: 11} \]

\noindent and we  say that $f$ is $\ell_{1}$-differentiable if the above limit exists.

\noindent \textbf{Theorem 2.2.}

\begin{enumerate} [label=(\alph*)]

\item  Let $f_k$ be $\ell_{1}$-differentiable for $k=0,1,2 ,... , n$, then

\[ ^{\ell_{1}} D(\sum\limits_{k=0}^{n} c_k f_k )(t) = \sum\limits_{k=0}^{n} c_k  (^{\ell_{1}} D f_k )(t). \tag{13} \label{eq: 12} \] 

\item Suppose that $f$ is $\ell_{1}$-differentiable, then

\[  ^{\ell_{1}} D(^{\ell_{1}} D f(t) ) =   {\lim\limits_{\epsilon \to 0}}  \mathscr{L}^{-1}(s^2 \mathscr{L}(|f(t-\epsilon)| - |f(\epsilon)| )) . \tag{12} \label{eq: 13} \]

\end{enumerate}

\noindent \textbf{Proof.} (a) is immediate, since

\[ ^{\ell_{1}} D( \sum\limits_{k=0}^{n} c_kf_k)(t) =  {\lim\limits_{\epsilon \to 0}} \mathscr{L}^{-1}(s \mathscr{L}(\sum\limits_{k=0}^{n} |c_k f_k(t-\epsilon)|) - |c_kf_k(\epsilon)| ) \] 

\[ = {\lim\limits_{\epsilon \to 0}} \mathscr{L}^{-1}(s\sum\limits_{k=0}^{n} \mathscr{L}( |c_k f_k(t-\epsilon)|) - |c_kf_k(\epsilon)| ) \] 

\[ = \sum\limits_{k=0}^{n} c_k {\lim\limits_{\epsilon \to 0}} \mathscr{L}^{-1}(s \mathscr{L}( | f_k(t-\epsilon)|) - |f_k(\epsilon)| ) \] 

\[ = \sum\limits_{k=0}^{n} c_k  (^{\ell_{1}} D ( f_k )). \]

\noindent To show (b), the result readily follows since

\[ = ^{\ell_{1}} D (  ^{\ell_{1}} D f)(t) =  {\lim\limits_{\epsilon \to 0}} \mathscr{L}^{-1}(s \mathscr{L}(| (^{\ell_{1}} D f)(t-\epsilon)|) - | (^{\ell_{1}} D f)(\epsilon)| )  \]

\[{\lim\limits_{\epsilon_1 \to 0}} \mathscr{L}^{-1}(s \mathscr{L}( \lim\limits_{\epsilon_2 \to 0} \mathscr{L}^{-1}(s  \mathscr{L} ( |f(t-\epsilon_1 - \epsilon_2)|) - |f(\epsilon_2-\epsilon_1)| )) - \lim\limits_{\epsilon_3 \to 0} \mathscr{L}^{-1} (s \mathscr{L} (|f(\epsilon_1-\epsilon_3)|) - |f(\epsilon_3)|)  \]

\[= {\lim\limits_{\epsilon_1 \to 0}} \mathscr{L}^{-1}(s \mathscr{L}( \lim\limits_{\epsilon_2 \to 0} \mathscr{L}^{-1}(s \mathscr{L} ( |f(t- \epsilon_1 - \epsilon_2)| )- |f(- \epsilon_1 - \epsilon_2)| )) - \mathscr{L}^{-1} (s \mathscr{L} (|f(\epsilon_1)| - |f(0)|) ) \]

\[ =  {\lim\limits_{\epsilon_1 \to 0}} \mathscr{L}^{-1}(s  \mathscr{L}( \lim\limits_{\epsilon_2 \to 0} \mathscr{L}^{-1}(s \mathscr{L} ( |f(t-\epsilon_1)|) - |f(\epsilon_1)| )) - \mathscr{L}^{-1} (s \mathscr{L} (|f(\epsilon_1)|) - |f(0)| ) \] 

\[ =  {\lim\limits_{\epsilon_1 \to 0}} \mathscr{L}^{-1}(s^{2} \mathscr{L} ( |f(t-\epsilon_1)|) - |f(\epsilon_1)| ) - \mathscr{L}^{-1} (s \mathscr{L} (|f(\epsilon_1)|) - |f(0)|)  \]

\[ = {\lim\limits_{\epsilon \to 0}}  \mathscr{L}^{-1}(s^2 \mathscr{L}(|f(t-\epsilon)| - |f(\epsilon)| )). \]

\noindent The above steps steps follow from Dominated Convergence, since

\[ |F(t)e^{-st}| \leq |F(t)|, \] 

\noindent and

\[ \int_{0}^{\infty} | F(t-\epsilon) - F(t) | dt  \to 0,  \]

\noindent which yields

\[  \lim\limits_{\epsilon \to 0} \mathscr{L}( |F(t-\epsilon)e^{-st}|) = \mathscr{L} (|F(t)e^{-st}|).  \] 

\bigskip

\qed

\noindent \textbf{Theorem 2.3.} Suppose that $f$ is real-valued, $\Phi: [0, \infty) \to \mathbb{R}$ is concave and $\ell_{1}$-differentiable with $\Phi(f(0))=0$. Then

\[ ^{\ell_{1}} (D)( \Phi(f) )  \geq \mathscr{L}^{-1} |\Phi(\mathscr{L}(f))|. \tag{14} \label{eq: 14}  \]  

\noindent \textbf{Proof.}

\noindent Jensen's inequality applied to $\Phi$ yields

 \[ \Phi(F(s-\epsilon)) \geq e^{-\epsilon} s^{-1}  \int_{0}^{\infty} se^{-st} \Phi(f(t-\epsilon)) dt, \] 
 
\noindent since 

\[ \int_{0}^{\infty} se^{-st} dt = 1. \]
 
\noindent The result follows from the definition of the $\ell_{1}$ derivative

\[ ^{\ell_{1}} (D)( \Phi(f) ) = {\lim\limits_{\epsilon \to 0}} \mathscr{L}^{-1}(s \mathscr{L}( |\Phi(f(t)-\epsilon) |) - s|\Phi(f(\epsilon)| ). \] 
 
 \qed
 
\noindent \textbf{Theorem 2.4.} Suppose $m \in \mathbb{R}$, $n \in \mathbb{N}$ and let $^{\ell_{1}} (D)(S)$ denote the space of $\ell_{1}$ differentiable functions on the set $S$. We have the following results:

\begin{enumerate} [label=(\alph*)]

\item $^{\ell_{1}} (D)(P_n) \subset P_{n-1}$, where $P_n(t)$ is the space of polynomials in $t$ of degree $n$.

\item Let $C \in \mathbb{R}$. Then 

\[  ^{\ell_{1}} (D)(C) =  |C|(\delta(t) - \delta'(t)). \tag{15} \label{eq: 15}  \] 

\item Let $m \in \mathbb{R}$. Then \[  ^{\ell_{1}} (D)(e^{mt})  =  \delta(t) - \delta'(t) + me^{m t}. \tag{16} \label{eq: 16} \]  

\item If $f$ is periodic with period $P>0$, differential and nonvanishing, then

 \[  ^{\ell_{1}} (D)(f(t))  = \frac{ f(t) f(t)' }{|f(t)|} - f(0). \tag{17} \label{eq: 17} \]

\item If $a >0$ and $f$ is of exponential order, differentiable and nonvanishing, then

\[ ^{\ell_{1}} (D)( f(at)  ) =  af'(at)|f(at)|^{-1}   - f(0)\delta'(t). \]

\end{enumerate} 

\noindent \textbf{Proof.} Let $k>0$. Property (a) follows from the computation

\[    ^{\ell_{1}} (D)(t^{2k})  =  {\lim\limits_{\epsilon \to 0}}  \mathscr{L}^{-1}(s \mathscr{L}( (t-\epsilon)^{2k} ) -s\epsilon^{2k} )    \]   

\[    =  {\lim\limits_{\epsilon \to 0}}  \mathscr{L}^{-1}( e^{\epsilon s} \Gamma(2k +1)  - s\epsilon^{2k})    \]   

\[    =  {\lim\limits_{\epsilon \to 0}} \mathscr{L}^{-1}(  e^{-\epsilon s} s^{-2k} \Gamma(2k +1)  - s\epsilon^{2k})    \]

\[    =  t^{2k - 1}    \]

\noindent To show (b),

\[ ^{\ell_{1}} (D)(C)=  {\lim\limits_{\epsilon \to 0}} \mathscr{L}^{-1}(s \mathscr{L}(|C|) - |C| )) \]

\[ = \mathscr{L}^{-1}(|C| - |C|s ) \]

\[ = |C|(\delta(t) - \delta'(t)),  \]

\noindent since $\mathscr{L} ^{-1} (s^{n})=D^{n} \delta(t)$ in the distributional sense.

\noindent Property (c) follows from the identity $\mathscr{L}( f(at) ) = \frac{F(s/a)}{a}$.

\[    ^{\ell_{1}} (D)(e^{mt})  =  {\lim\limits_{\epsilon \to 0}}   \mathscr{L}^{-1}(s \mathscr{L}( e^{m(t-\epsilon)}) - se^{\epsilon})    \]  

\[    =  {\lim\limits_{\epsilon \to 0}}   \mathscr{L}^{-1}(sm^{-1} (sm^{-1} -1)^{-1}e^{(\epsilon sm^{-1} )}  - se^{m\epsilon})    \] 

\[    =  {\lim\limits_{\epsilon \to 0}}   \mathscr{L}^{-1}(s(s-1)^{-1}e^{\epsilon s}  - se^{m\epsilon})    \]  

\[    = \delta(t) - \delta'(t) + me^{m t}.  \]

\noindent To show (d), applying the identity

\[ \mathscr{L} (|f(t-\epsilon)|) = \frac{e^{s\epsilon} \mathscr{L}(f)}{1-e^{-sP}}, \]

\noindent we obtain

\[ ^{\ell_{1}} (D)(f(t)) = \lim\limits_{\epsilon \to 0^{+} } \mathscr{L}^{-1} ( ( s^{1 + \frac{s \epsilon}{ln s} } ) (1-e^{-sP})^{-1} \mathscr{L}(f) - s|f(\epsilon)|).  \]

\noindent To show (e), we have

\[ ^{\ell_{1}} (D)( f(at)  ) =  {\lim\limits_{\epsilon \to 0}}  \mathscr{L}^{-1}(s \mathscr{L}( f(a(t-\epsilon)) - s|f(a\epsilon)| ) \]

\[ =  {\lim\limits_{\epsilon \to 0}}  \mathscr{L}^{-1}(sa^{-1} e^{a\epsilon s} \mathscr{L}(| f(at)|) - s|f(a\epsilon)| ) \]

\[ =  {\lim\limits_{\epsilon \to 0}}  \mathscr{L}^{-1}(sa^{-1} e^{a\epsilon s} \mathscr{L}(F^{*}(s/a)) - |f(a\epsilon)| ) \]

\[  = af'(at)|f(at)|^{-1}   - f'(0)\delta'(t).   \]

 \qed

\bigskip

\bigskip

\noindent We now consider generalizations of the work above. One consideration is to adjust the rate of convergence in ($\ref{eq: 11}$).  Let $f$ be $\ell_{1}$-differentiable. Define the modified $\ell_{1}$-fractional derivative of order $\alpha \in \mathbb{R}$ by

\[  ^{\ell_{1} +}D^{\alpha}_{t} f = \lim\limits_{\epsilon \to 0^{+}} \mathscr{L}^{-1} \bigg( s^{\alpha} \big( \mathscr{L}(f(t-\epsilon) - f(t\epsilon)) \bigg).   \tag{19} \label{eq: 19} \] 

\noindent In the case that $f(t)=|t|$, we obtain

\[ \mathscr{L}^{-1}(\epsilon s^{\alpha - 1} +  (1-\epsilon + 2e^{-\epsilon s})s^{\alpha-2} ). \]
\noindent When $\alpha = 1$, this resembles the $\ell_{1}$ fractional derivative. A quick application of the convolution theorem, leads to a second representation

\[  ^{\ell_{1} +}D^{\alpha}_{t} f = \lim\limits_{\epsilon \to 0^{+}} \mathscr{L}^{-1} \bigg( s^{\alpha} \big( \mathscr{L}(f(t-\epsilon) - f(t\epsilon) \bigg).   \tag{19} \label{eq: 19}  =  \lim\limits_{\epsilon \to 0^{+}} \mathscr{L}^{-1} \bigg(s^{\alpha - 1} s \bigg( \mathscr{L}\bigg( f(t - \epsilon) - f(s \epsilon) \bigg) \bigg)   \] 

\[ =  \lim\limits_{\epsilon \to 0^{+}} \frac{t^{-\alpha}}{\Gamma(1-\alpha)}  \star (f(t - \epsilon) - f(t \epsilon))'  =  \frac{1}{\Gamma(1-\alpha)} \lim\limits_{\epsilon \to 0^{+}}  \int_{0}^{t} (s-t)^{-\alpha} \bigg( f'(s-\epsilon) - \epsilon f'(s\epsilon) \bigg) \hspace{1mm} ds. \]

\noindent A natural generalization of the derivative above can be seen by considering the generalized $\beta,\gamma$ finite difference operator defined for $0<a<t<b$ , $f \in H^{1}(a,b)$, $\alpha \in \mathbb{R}$ by

\[  ^{\beta, \gamma}D^{\alpha}_{t} f  =  \frac{1}{Mf(1-\alpha)} \lim\limits_{\epsilon \to 0^{+}}  \int_{a}^{t} (s-t)^{-\alpha} \sum\limits_{k=1}^{\ceil{\alpha}} \bigg( f'(s^{k}-\epsilon) - kE_{\beta,\gamma}(\log \epsilon) f'(s\epsilon^{k}) \bigg) \hspace{1mm} ds, \] 

\noindent where $Mf(s)$ denotes the Mellin transform of $f$ and  $E_{\beta, \gamma}$ denotes the two-parameter Mittag-Leffler function given by

\[ E_{\beta, \gamma} (z) = \sum\limits_{k=0}^{\infty} \frac{z^{k}}{\Gamma(\beta k + \gamma )}, \tag{20} \label{eq: 20} \] 
 
\noindent  where $\beta, \gamma \in \mathbb{C}$. In Oliveira \cite{Ref23}, the following generalized Laplace fractional derivative is discussed. Suppose that $\Phi(s,\alpha)$ is a two-parameter real-valued operator satisfying

\[ \mathscr{L} (_{0} D^{\alpha}_{t} f) = \Phi(s,\alpha) \mathscr{L}(f(t)) (s), \tag{21} \label{eq: 21} \]

\noindent where $D^{\alpha}$ denotes a general fractional derivative. If $\Phi(s,\pm 1)  = s^{ \pm 1}$ and $\Phi(s,0)=1$ and the operator $\Phi$ also satisfies $\Phi(s, \alpha)=sk(s,\alpha)$, proceeding as in ($\ref{eq: 11}$), we are left with 

\[_{0} D^{\alpha}_{t} f = \int_{0}^{t} K(t-x, \alpha) Df(x) dx, \tag{22} \label{eq: 22} \] 

\noindent where $K(t,\alpha) = \mathscr{L}^{-1} (k(s,\alpha))$. In ($\ref{eq: 22}$), the kernel $K(t,\alpha)$ plays the role of $g_{n-\alpha}$ in ($\ref{eq: 11}$). Caputo and Fabrizio \cite{Ref22} propose a Caputo-type fractional derivative with an exponential kernel, considering applications to constitutive equations for dissipation in \cite{Ref10}. This derivative is studied in distributional settings in Atanakovic \cite{Ref26}, where this operator is shown to obey a viscoelastic consistency result. The Caputo-Fabrizio fractional derivative is criticized in Tarasov \cite{Ref11} and Origuieira et al. \cite{Ref12} for its non-locality. In Ortiguieira et al. \cite{Ref30}, criteria for suitable fractional derivative candidates are considered including reduction to a classical derivative, backwards compatibility, a neutral element and the Leibniz rule.
\bigskip

A generalized Laplace transform approach is considered in Oliveira \cite{Ref23} working from the theory of distributions, revealing a deep relationship between the Riemann-Liouville and Caputo fractional derivative definitions. This leads the authors back to the choice of kernel in \cite{Ref22}, which is then viewed in the context of a relaxation model. The interaction between the Laplace transform and Mittag-Leffler functions can be found in Gorenflo and Mainardi \cite{Ref9}.  A generalized fractional derivative is also seen in Katugampola \cite{Ref27}, extending the limit definition of the derivative. It is shown to satisfy linearity, the Leibniz rule and the Mean Value Theorem. The authors also introduce a generalized fractional integral satisfying the integration by parts formula. In \cite{Ref28}, a generalization of the Riemann-Liouville and Hadamard fractional derivatives is investigated.

\bigskip

In the context of gradient descent algorithms, it is possible to replace the usual gradient $\nabla f$ with a fractional gradient based on the modified $\ell_{1}$ fractional derivative. We have seen that even for elementary functions, the $\ell_{1}$ fractional derivative often exists only in the distributional sense, which means that convergence to an extreme point is more difficult to achieve. In the context of $\ell_{1}$ regularization, it is evident that this fractional gradient mimics the classical gradient however, this may not be true of a generalization. The search for a fractional gradient with better smoothing properties is related to the choice of kernel in ($\ref{eq: 11}$).

\section{Conclusion}

In this paper, we have discussed a new derivative defined by the Laplace transform. This representation interacts well with periodic functions and satisfies some of the expected classical properties. This derivative can be readily generalized to one of fractional order. We have investigated a more suitable operator with better smoothing properties that are compatible with $\ell_{1}$ norm problems. The generalization in ($\ref{eq: 21}$) and the viewpoint of convolution in ($\ref{eq: 11}$) is useful in the search for a better kernel. As a future work, we seek a fractional operator that is nonlocal yet manages to smoothen the $\ell_{1}$ loss function in gradient descent methods.

\bibliographystyle{unsrt}  






\end{document}